\documentclass{article}
\usepackage{graphicx} 
\usepackage{subcaption} 
\usepackage{color}
\usepackage[colorlinks=true,allcolors=blue
]{hyperref}
\usepackage{amsmath, amssymb, amsthm}
\usepackage{mathrsfs}
\usepackage{mathtools}
\usepackage{fullpage}
\usepackage[noadjust,nocompress]{cite}
\usepackage{graphics}
\usepackage{xcolor}
\usepackage{tikz}

\usepackage{thmtools}
\usepackage[noabbrev,capitalize,nameinlink]{cleveref}

\newtheorem{theorem}{Theorem}
\newtheorem{proposition}[theorem]{Proposition}

\newtheorem{lem}[theorem]{Lemma}

\newtheorem{claim}[theorem]{Claim}
\newtheorem{corollary}[theorem]{Corollary}

\newtheorem{ques}[theorem]{Question}

\newcommand*{\myproofname}{Proof}
\newenvironment{claimproof}[1][\myproofname]{\begin{proof}[#1]}{\end{proof}}

\newcommand*{\abs}[1]{\lvert #1\rvert}

\title{A general bound on $R(C_k,H)$} 
\author{Stijn Cambie\thanks{Department of Computer Science, KU Leuven Campus Kulak-Kortrijk, 8500 Kortrijk, Belgium. Supported by a postdoctoral fellowship by the Research Foundation Flanders (FWO) with grant number 1225224N. Email: \protect\href{mailto:stijn.cambie@hotmail.com}{\protect\nolinkurl{stijn.cambie@hotmail.com}}} \and Andrea Freschi
\thanks{HUN-REN, Alfr{\'e}d R{\'e}nyi Institute of Mathematics, Budapest, Hungary. Research partially supported by ERC Advanced Grants ``GeoScape", no. 882971 and ``ERMiD", no. 101054936. E-mail: \protect\href{freschi.andrea@renyi.hu}{\protect\nolinkurl{freschi.andrea@renyi.hu}}} 
}

\date{}

\begin{document}

\maketitle
\begin{abstract}
    In this paper, we prove that for every $k$ and every graph $H$ with $m$ edges and no isolated vertices, the Ramsey number $R(C_k,H)$ is at most $(k-1)m+1\le km$.  
    This settles a problem of Erd{\H{o}}s, Faudree, Rousseau and Schelp, which is listed as problem 34 in the graph theory collection.
\end{abstract}

\section{Introduction}
Ramsey theory is a central research topic in combinatorics.
In his seminal work, Ramsey \cite{ramsey1930} proved that for any pair of graphs~$G$ and~$H$ there exists a positive integer~$N$ such that any red-blue edge-colouring of the complete $N$-vertex graph contains a red~$G$ or a blue~$H$. 
The {\it Ramsey number}~$R(G,H)$ is the smallest~$N$ for which the above holds. 

Let~$C_k$ denote the cycle of length~$k$.
Goddard and Kleitman, and independently Sidorenko, proved the following result.

\begin{theorem}[Goddard and Kleitman~\cite{goddard1994upper}, and Sidorenko~\cite{Sidorenko93}]\label{thm:C_3}
    For every graph~$H$ with~$m$ edges and no isolated vertices, we have $R(C_3,H)\le 2m+1$.
\end{theorem}
The bound in \cref{thm:C_3} is best possible, as equality is attained when~$H$ is a tree or a matching.

\smallskip

The study of~$R(C_k,H)$ for other values of~$k$ was later considered in a number of works, we refer to the introduction in~\cite{CFMPP26} for a more detailed history of these problems.
In~$1993$, Erd{\H{o}}s, Faudree, Rousseau and Schelp raised the following question.

\begin{ques}[Erd{\H{o}}s, Faudree, Rousseau and Schelp~\cite{EFRS93}]\label{question}
    Given an odd integer $k\ge3$, what is the smallest~$c_k$ such that $R(C_k,H)\le c_km$ for every graph~$H$ with~$m$ edges an no isolated vertices? 
\end{ques}

This problem was also popularised as problem $34$ in the book by Chung and Graham~\cite{chung1998erdos} and in the database of difficult graph theory Erd\H{o}s problems~\cite{ChungErdosProblems} from 2010, and as problem $569$ in the (intended to be) exhaustive database of Erd\H{o}s problems~\cite{BloomErdosProblems569}.

Since $R(C_k, K_2)=k,$ one immediately notices $c_k \ge k$.
In this article, we prove that~$c_k=k$ is the correct constant.
In fact, we obtain a slightly stronger general upper bound.

\begin{theorem}\label{thm:main}
    For every graph $H$ with $m$ edges and no isolated vertices, and for every $k \ge 3,$
    $$R(C_k,H) \le (k-1) m+1 \le km.$$
\end{theorem}

Note that \cref{thm:C_3} corresponds to the case~$k=3$ of \cref{thm:main}.
Therefore, the bound~$(k-1)m+1$ is tight for arbitrary $k$ and $H=K_2$, as well as when $k=3$ and~$H$ is a tree or a matching.

\smallskip

Observe that sharper bounds on~$R(C_k,H)$ are known to hold under additional assumptions.
Recently~\cite{CFMPP26}, it was proved that $R(C_k,H)\le 2m+\lfloor(k-1)/2\rfloor$ provided~$m$ is sufficiently large with respect to~$k$.
On the other hand, Burr~\cite{Burr} proved that $R(C_k,H)=(\chi(H)-1)(k-1)+\sigma(H)$,
\footnote{Here~$\chi(H)$ and~$\sigma(H)$ denote the chromatic number and the smallest colour size in a proper~$\chi(H)$-colouring of~$H$, respectively.}
provided~$k$ is large with respect to~$H$ (for quantitative bounds on~$k$ see e.g.,~\cite{PS20,HHKL23}).
The bound in \cref{thm:main} is thus not optimal in general, but it has no restriction on~$k$ and~$m$. 

\smallskip

To prove \cref{thm:main}, we follow the same ideas as in~\cite{CFMPP26}, with some appropriate modifications.
The proof of \cref{thm:main} is discussed in \cref{sec:mainproof}. 
In \cref{sec:end} we discuss related questions and some observations that arose while working on this note.

\subsection{Notation}

We use standard terminology and notation in graph theory.
We denote the path on~$k$ vertices as~$P_k$.
For a graph $G=(V,E)$ and subset $U \subseteq V$ we let $G[U]$ denote the graph induced by $U$ i.e., $G[U]:=(U, E \cap \binom{V}{2})$.
For a vertex $v\in V(G)$, we write $G-v$ to denote the graph obtained by removing~$v$ from~$G$.
The chromatic number of~$G$ is denoted by $\chi(G).$
We write~$e(G)$ and~$|G|$ to denote the number of edges and vertices of~$G$, respectively.
We write~$\delta(G)$ and~$\bar d(G)$ for the minimum and average degree of~$G$, respectively.
The second neighbourhood of a vertex $v\in V(G)$ is the set of vertices in~$G$ at distance exactly~$2$ from~$v$.
For a positive integer~$n$, we write $[n]:=\{1,2,\dots,n\}$ for the set of the first~$n$ positive integers.

\section{Proof of~\cref{thm:main}}\label{sec:mainproof}

We first list some preliminary results.
The next lemma addresses some instances where~$k$ is small, and follows easily by combining several results in the literature.


\begin{lem}\label{lem:jayawardene}
    For every connected graph $H$ with $m$ edges and no isolated vertices, and for every $k \in\{4,5,6\}$, we have
    $$R(C_k,H) \le (k-1) m+1.$$
\end{lem}
\begin{proof}
    The statement of the lemma holds for $m=1$, so we may assume $|H|\ge3$.
    By the work of Jayawardene~\cite[Thms.~4.1,~4.5,~4.7]{jayawardene1999thesis} we have
        $$R(C_k,H)=\begin{cases}
        2m+1 \text{ for $k=4$ and $m\ge 2$},\\
        2m+2 \text{ for $k=5$ and $|H|\ge 4$},\\
        2m+2 \text{ for $k=6$ and $|H|\ge 6$}.
    \end{cases}$$
    Hence, it remains to verify the cases $k=5$ and $|H|=3$, and $k=6$ and $3\le|H|\le 5$.
    These follow from $R(C_5,K_3)=9$, $R(C_6,K_3)=11$, $R(C_6,K_4)=16$ and $R(C_6,K_5)=21$ (see the dynamic survey by Radziszowski~\cite{Radziszowski}) and the fact that $m\ge n-1$ (since~$H$ is connected).
\end{proof}

Next, we state a sharper version of~\cite[Lem.~6]{CFMPP26}.
The proof, very similar to that in~\cite{CFMPP26}, can be found in~\cref{app:proof_lem3}.

\begin{lem}\label{lemma:P_kvsH}
    For every integer $k\ge 2$ and graph~$H$, we have $R(P_k,H)\le |H|+(k-2)(\chi(H)-1)$.
\end{lem}

Using $\chi(H)\le |H|$, we immediately deduce the following corollary.

\begin{corollary}\label{col:P_kvsH}
    For every integer $k\ge 2$ and graph~$H$, we have $R(P_k,H)\le (k-1)(|H|-1)+1$.
\end{corollary}

We also use the following strengthening of \cite[Lem.~8]{CFMPP26}.

\begin{lem}\label{lemma:path_second_neighbourhood_stronger}
    Let $k \ge 7$. For any graph $G$ and vertex $v\in V(G)$, if the second neighbourhood of $v$ contains a copy of $P_{k+1}$ then $G$ contains a copy of $C_k$.
\end{lem}

\begin{proof}
    Let the vertices of the copy of $P_{k+1}$ be, in order, $v_1, v_2, \ldots, v_{k}, v_{k+1}$.
    For every $i\in [k+1]$, let $u_i$ be an arbitrary neighbour of $v$ adjacent to $v_i$.

    \begin{claim}\label{claim:k-4+i}
        For every $i\in[5]$, we have $u_i=u_{k-4+i}$.    
    \end{claim}

    \begin{claimproof}
        If $u_i\neq u_{k-4+i}$ then $vu_iv_i\dots v_{k-4+i}u_{k-4+i}v$ is a copy of~$C_k$, contradiction.
    \end{claimproof}

    \begin{claim}\label{claim:+}
        For every $i\in[3]$, we have $u_i\neq u_{i+2}$.   
    \end{claim}

    \begin{claimproof}
        By Claim~\ref{claim:k-4+i}, we have $u_{i+2}=u_{k-2+i}$.
        
        If $u_i=u_{i+2}=u_{k-2+i}$ then $u_iv_i\dots v_{k-2+i}u_{k-2+i}$ is a copy of~$C_k$, see~\cref{fig1:a}.
        Thus we have $u_i\neq u_{i+2}$. \end{claimproof}

Now we are ready to finish the proof, taking into account the above lemmas, by considering a few cases.

    If $k=7$, we have $u_1\neq u_3$, $u_3=u_6$ and $u_1=u_4=u_7$ and so $u_1v_1v_2v_3u_3v_6v_7u_1$ is a copy of $C_k.$

    If $k \ge 8$ and $u_1=u_{5}$, then $u_1v_1v_{2}v_{3}u_{3}v_{k-1}v_{k-2}\dots v_{5}u_{1}$ is a copy of~$C_k$, see~\cref{fig1:b}.

    If $k \ge 8$ and $u_1 \not= u_5,$ then $u_1v_1v_{2}v_{3}u_{3}v u_5 v_5v_6 \ldots v_{k-3} u_1$ is a copy of~$C_k$, see~\cref{fig1:c}.

    \begin{figure}[h!]   
    \centering
        \begin{subfigure}{.32\textwidth}
        \centering
        \begin{tikzpicture}[scale=0.94]
        \draw (1.5,3.5) -- (4,0) -- (-1,0) -- (1.5,3.5);
        \filldraw (1.5,3.5) circle (2pt) node[anchor=south]{$v$};
        \filldraw[fill=white] (1.5,2) ellipse (1.5cm and 0.7cm);
        \filldraw[fill=white] (1.5,0) ellipse (2.5cm and 1cm);
        \filldraw[very thick] (0,0) circle (2pt) node[anchor=north]{$v_i$} -- (1,0) circle (2pt) node[anchor=north]{$v_{i+1}$} -- (1.5,0);
        \filldraw[very thick,dotted] (1.5,0) -- (2.5,0);
        \filldraw[very thick] (3,0) circle (2pt) node[anchor=north]{$v_{k-2+i}$} -- (2.5,0);
        \filldraw[very thick] (0,0) -- (1.5,2) circle (2pt) node[anchor=south]{$u_i=u_{k-2+i}$} -- (3,0);   
        \end{tikzpicture}
        \caption{}
        \label{fig1:a}
        \end{subfigure}   
        \hfill
        \begin{subfigure}{.32\textwidth}
        \centering
        \begin{tikzpicture}[scale=0.94]
        \draw (1.5,3.5) -- (4,0) -- (-1,0) -- (1.5,3.5);
        \filldraw (1.5,3.5) circle (2pt) node[anchor=south]{$v$};
        \filldraw[fill=white] (1.5,2) ellipse (1.5cm and 0.7cm);
        \filldraw[fill=white] (1.5,0) ellipse (2.5cm and 1cm);
        \filldraw[very thick] (1,2) circle (2pt) node[anchor=east]{$u_{1}$} -- (-0.2,0) circle (2pt) node[anchor=north]{$v_1$} -- (0.3,0) circle (2pt) -- (0.8,0) circle (2pt) node[anchor=north]{$v_{3}$} -- (2,2) circle (2pt) node[anchor=west]{$u_{3}$} -- (3.3,0) circle (2pt) node[anchor=north]{$v_{k-1}$} -- (2.9,0);
        \filldraw (1.3,0) circle (2pt);
        \filldraw[very thick,dotted] (2.2,0) -- (2.9,0);
        \filldraw[very thick] (2.2,0) -- (1.8,0) circle (2pt) node[anchor=north]{$v_{5}$} -- (1,2);  
        \end{tikzpicture}
        \caption{}
        \label{fig1:b}
        \end{subfigure}
        \hfill
        \begin{subfigure}{.32\textwidth}
        \centering
        \begin{tikzpicture}[scale=0.94]
        \draw (1.5,3.5) -- (4,0) -- (-1,0) -- (1.5,3.5);
        \filldraw (1.5,3.5) circle (2pt) node[anchor=south]{$v$};
        \filldraw[fill=white] (1.5,2) ellipse (1.5cm and 0.7cm);
        \filldraw[fill=white] (1.5,0) ellipse (2.5cm and 1cm);
        \filldraw[very thick] (1.5,3.5)--(1,2) circle (2pt) node[anchor=east]{$u_3$} -- (0.8,0) circle (2pt) node[anchor=north]{$v_3$} -- (0.3,0) circle (2pt) -- (-0.2,0) circle (2pt) node[anchor=north]{$v_{1}$} -- (1.5,2) circle (2pt) node[anchor=south]{$u_{1}$} -- (3.3,0) circle (2pt) node[anchor=north]{$v_{k-3}$} -- (2.9,0);
        \filldraw (1.3,0) circle (2pt);
        \filldraw[very thick,dotted] (2.2,0) -- (2.9,0);
        \filldraw[very thick] (2.2,0) -- (1.8,0) circle (2pt) node[anchor=north]{$v_{5}$} -- (2,2) circle (2pt) node[anchor=west]{$u_{5}$}-- (1.5,3.5);  
        \end{tikzpicture}
        \caption{}
        \label{fig1:c}
        \end{subfigure} 
        \vspace*{-0.25cm}
        \caption{From left to right, the cases $u_i=u_{i+2}$, $u_1=u_{5}$ and $u_1\not=u_{5}$ in the proof of~\cref{lemma:path_second_neighbourhood_stronger}.}
        \label{fig1}
    \end{figure}
\end{proof}

Finally, we use the following elementary relation between the number of edges, vertices and the chromatic number of a graph.

\begin{proposition}\label{prop:lowbnd_m}
    Let $G$ be a graph with~$n$ vertices, $m$ edges and chromatic number~$\chi$.
    If~$\delta(G)\ge2$ then $m \ge n+ \frac{(\chi-3) \chi}{2}$.
\end{proposition}

\begin{proof}
    Since $\delta(G)\ge2$, we have $m\ge n$.
    This implies the statement of the proposition for~$\chi=2$, thus we may assume that $\chi\ge 3$.
    
    Let~$H$ be a {\it vertex-colour-critical} subgraph of~$G$, that is, such that~$\chi(H)=\chi$ and~$\chi(H')<\chi$ for every proper subgraph~$H'$ of~$H$.
    Note that $|H|\ge\chi$.
    Furthermore, if a vertex~$v\in V(H)$ has degree at most~$\chi-2$ in~$H$, then we could extend a proper $(\chi-1)$-colouring of~$H-v$ to~$H$, a contradiction.
    Hence, $\delta(H)\ge\chi-1$.
    The bounds $\delta(H)\ge\chi-1$ and $\delta(G)\ge 2$ imply
    $$\sum_{v\in V(G)}d(v)=\sum_{v\in V(H)}d(v)+\sum_{v\notin V(H)}d(v)\ge|H|(\chi-1)+2(n-|H|)=|H|(\chi-3)+2n.$$
    By the hand-shaking lemma and $|H|\ge\chi\ge 3$, it follows that $2m\ge\chi(\chi-3)+2n$.
    Dividing by~$2$ yields the claimed inequality.
\end{proof}

Now we are ready to prove our main theorem.


\begin{proof}[Proof of \cref{thm:main}]
    We proceed by induction on~$e(H)$ with~$k$ fixed. 
    For the base case, if~$e(H)=1$ then~$H$ consists of a single edge and so $R(C_k,H)=k$, as required.
    Thus, we may assume~$e(H)>1$.
    Let $n:=|H|$ and $m:=e(H)$.
    For a contradiction, suppose there exists a red-blue edge-coloured clique~$G$ on~$N:=(k-1)m+1$ vertices which does not contain a red~$C_k$ nor a blue~$H$.
    
    If~$H$ is disconnected then~$H$ is the disjoint union of two graphs~$H_1$ and~$H_2$, each with no isolated vertices.
    By inductive hypothesis we have $R(C_k,H_1)\le (k-1)\cdot e(H_1)+1$, and so~$G$ contains a blue copy of $H_1$. Since $N-|H_1| \ge N-2e(H_1)\ge(k-1)\cdot e(H_2)+1$, then $G$ contains a blue copy of~$H_2$ which is vertex-disjoint from the blue copy of~$H_1$.
    Thus~$G$ contains a blue copy of~$H$, contradiction.

    Therefore, we may assume~$H$ is connected, and so $m\ge n-1$.
    If $k=3$ or $4\le k\le6$ then~$G$ contains a red~$C_k$ or a blue~$H$ by \cref{thm:C_3} and \cref{lem:jayawardene} respectively, a contradiction.
    Thus, we may assume $k\ge 7$.
    Also, observe that $e(H)>1$ implies $n\ge 3$.

    \smallskip

    Let~$u\in H$ be a vertex of minimum degree in~$H$.
    Since~$H$ is connected and $e(H)>1$, it follows that~$H-u$ does not have isolated vertices and~$e(H-u)<e(H)$.
    By inductive hypothesis, $G$ contains a blue copy~$F$ of~$H-u$.
    If there is a vertex outside $F$ that is connected with only blue edges to the images of the neighbours of~$u$ then~$G$ contains a blue copy of~$H$, a contradiction. 
    So we conclude there is a red star, say with centre~$v$, and a set~$S$ of leaves with $|S|\ge (N-n+1)/\delta(H)$.

    The subgraph~$G[S]$ does not contain a red $P_{k-1}$, as that would form a red~$C_k$ with~$v$.
    Also, $G[S]$ does not contain a blue $H$.
    If $\delta(H)=1$, we have 
    $$|S|\ge (k-1)m+1-n+1\ge(k-2)(n-1)+1$$
    where we used that $m\ge n-1$.
    By \cref{col:P_kvsH}, $G[S]$ contains a red~$P_{k-1}$ or a blue~$H$, contradiction. 
    Hence, we may assume $\delta(H)\ge2$. 

    \smallskip
    
    From the hand-shaking lemma we have $2m=n\cdot\bar d(H)$.
    It follows that
    \begin{equation}\label{eq:sizeS}
        \abs S \ge \frac{(k-1)m+1-n+1}{\delta(H)}\ge\frac{(k-1)m-n+2}{\bar d(H)}=\frac{(k-1)n}{2}-\frac{n-2}{\bar d(H)}\ge \frac{(k-2)n}{2}+1
    \end{equation}
    where the last inequality follows from $\bar d(H)\ge\delta(H)\ge2$.
    It follows from Corollary~\ref{col:P_kvsH} that~$G[S]$ contains a red~$P_{k-1}$ or a blue clique on~$\lfloor n/2\rfloor+1$ vertices.
    Note that we excluded the former already.
    Thus, let~$U_1\subseteq S$ with~$|U_1|=\lfloor n/2\rfloor+1$ and so that~$G[U_1]$ is a blue clique.

    Let~$\Gamma$ be the set of vertices in~$V(G)\setminus(\{v\}\cup U_1)$ incident to some red edge which is also incident to~$U_1$.
    Clearly, $G[\Gamma]$ does not contain a blue~$H$.
    Furthermore, $G[\Gamma]$ does not contain a red~$P_{k+1}$ by Lemma~\ref{lemma:path_second_neighbourhood_stronger} since~$\Gamma$ is the second neighbourhood of~$v$ within the graph spanned by the red edges between~$v$ and~$U_1$ and the red edges between~$U_1$ and~$\Gamma$.
    Hence, by \cref{lemma:P_kvsH}
    $$|\Gamma|\le R(P_{k+1},H)-1\le (n-1)+(k-1)(\chi(H)-1).$$

    Let~$U_2:=V(G)\setminus(U_1\cup\Gamma\cup\{v\})$.
    Then all edges between~$U_1$ and~$U_2$ are blue, by the definition of~$\Gamma$.
    Take a bipartition $H_1\sqcup H_2=V(H)$ uniformly at random with $|H_1|=|U_1|=\lfloor n/2\rfloor+1$ and $|H_2|=\lceil n/2\rceil-1$.
    Note that $\lceil n/2\rceil-1\ge1$ since $n\ge 3$.
    We have
    \begin{align*}
    \mathbb E(e(H_2))=m\cdot\frac{\lceil n/2\rceil-1}{n}\cdot\frac{\lceil n/2\rceil-2}{n-1}&\le m\cdot\frac{n-1}{2n}\cdot\frac{n-3}{2(n-1)}=\frac{m}{4}\left(1-\frac{1}{n}\right)\left(1-\frac{2}{n-1}\right)<\frac{m}{4}.
    \end{align*}
    
    Fix a bipartition $H_1\sqcup H_2=V(H)$ such that $|H_1|=\lfloor n/2\rfloor+1$ and $|H_2|=\lceil n/2\rceil-1$ and~$e(H_2)< m/4$.
    Observe that~$G[U_2]$ does not contain a red~$C_k$.
    Also, it does not contain a blue~$H_2$, since that would imply $G[U_1\cup U_2]$ contains a blue~$H$.
    Thus, we have $|U_2|\le R(C_k,H_2)-1$.
    By inductive hypothesis, we have $R(C_k,H_2)\le \max\{(k-1)\cdot e(H_2)+1,|H_2|\}$.
    Note that $|H_2|\le n/2\le m/2\le m(k-1)/4$.
    Hence,
    $$|U_2|<\frac{m(k-1)}{4}.$$   
    It follows that
    \begin{equation*}
        (k-1)m+1=1+|U_1|+|U_2|+|\Gamma|< 1+\left(\frac{n}{2}+1\right)+\left(\frac{m(k-1)}{4}\right)+\left((n-1)+(k-1)(\chi(H)-1)\right).
    \end{equation*}
    Rearranging yields
    \begin{equation*}
        \frac{3(k-1)m}{4}< \frac{3n}{2}+(k-1)(\chi(H)-1).
    \end{equation*}
    Since $k\ge 7$, we have
    \begin{equation}\label{eq:final_1}
        \frac{3m}{4}< \frac{3n}{2(k-1)}+\chi(H)-1\le\frac{n}{4}+\chi(H)-1.
    \end{equation}
    On the other hand, by \cref{prop:lowbnd_m} we have
    \begin{equation}\label{eq:final_2}
        \frac{3m}{4}\ge \frac{3n}{4}+\frac{3}{8}\chi(H)(\chi(H)-3).
    \end{equation}
    Combining~\eqref{eq:final_1} and~\eqref{eq:final_2} yields
    $$\chi(H)-1>\frac{n}{2}+\frac{3}{8}\chi(H)(\chi(H)-3).$$
    For $n\ge4$, the above inequality fails for any integer value of $\chi(H)$.
    Hence we have $n=3$.
    Since $\delta(H)\ge2$, it follows $H=C_3$.
    By \cref{thm:C_3}, it follows $R(C_k,C_3)\le 2k+1\le N$, and so $G$ contains either a red $C_k$ or a blue $C_3$.
    This is a contradiction, and concldues the proof.
\end{proof}

\section{Concluding remarks}\label{sec:end}

There are several interesting open questions about the Ramsey number~$R(C_k,H)$.
For example, Erd{\H{o}}s, Faudree, Rousseau and Schelp~\cite{EFRS78} conjectured that $R(C_k,K_n)=(n-1)(k-1)+1$ for all $k \ge n\ge3$ where $(n,k) \not= (3,3)$.
This is known to be true asymptotically, see~\cite{KLS21}. 
Erd\H{o}s \cite{ErdosP} posed the question of proving that $R(C_4,K_n)\le n^{2-\epsilon}$ for some $\epsilon>0$ and every $n\ge1$; this is still open.
Next to these two problems which are also listed as problems $17$ and $18$ in~\cite{ChungErdosProblems} (\#159 and \#551 in \url{erdosproblems.com}), also $\#19$ and $\#31, \#32, \#33$ are open and have similar flavour.

\smallskip

We remark that, at least for $k$ somewhat large, \cref{thm:main} can be proven using different ideas.
Recall the partition $H_1\sqcup H_2$ of~$V(H)$ in the proof of \cref{thm:main} was picked uniformly at random (subject to fixed sizes of~$H_1$ and~$H_2$).
Instead, one can select~$H_1$ to be an independent set incident to many edges.
Since~$H_1$ is independent, it can be easily embedded into~$S$.
On the other hand, one has a bound on the number of edges within~$H_2$, which allows us to invoke the inductive hypothesis (similarly to how we bounded~$|U_2|$).
See \cref{app:edge_heavy} for a precise statement about independent sets incident to many edges.


\appendix

\crefalias{section}{appendix}

\section{Proof of~\cref{lemma:P_kvsH}}\label{app:proof_lem3}

\begin{proof}
    It suffices to prove that for every integers $k\ge2$, $t\ge 1$ and $n_1,\dots,n_t\ge1$ we have    \begin{equation}\label{eq:RamseyPathMultipartite}
        R(P_k, K_{n_1, n_2, \ldots, n_t}) \le (k-2)(t-1)+ \sum_{i=1}^t n_i
    \end{equation}
    where~$K_{n_1, n_2, \ldots, n_t}$ denotes the complete $t$-partite graph with parts of size $n_1,\dots,n_t$.
    Indeed~\eqref{eq:RamseyPathMultipartite} implies the statement of the lemma, since~$H$ is a subgraph of $K_{n_1, n_2, \ldots, n_t}$ where $t= \chi(H)$ and the $n_i$ are the sizes of the colour classes in a proper $t$-colouring of~$V(H)$.

    We prove~\eqref{eq:RamseyPathMultipartite} by induction on~$t$. 
    The case $t=1$ is trivial.
    The case $t=2$ follows from a result of H\"aggkvist~\cite{Haggkvist89}, namely that $R(P_k,K_{n_1,n_2})=n_1+n_2+k-2$.
    Suppose~$t\ge3$.
    Given a blue copy of~$K_{n_1+n_2+k-2, n_3, \ldots, n_t}$ within a red-blue edge-coloured complete graph,
    we can apply H\"aggkvist's result to the graph spanned by the part of order~$n_1+n_2+k-2$ and find either a red copy of~$P_k$ or a blue copy of~$K_{n_1,n_2, n_3, \ldots, n_t}$.
    By this observation, and by the inductive hypothesis, we have
    \begin{align*}
    R(P_k, K_{n_1, n_2, \ldots, n_t}) \le R(P_k, K_{n_1+n_2+k-2, n_3, \ldots, n_t})&\le (k-2)(t-2)+(n_1+n_2+k-2)+\sum_{i=3}^t n_i\\
    &=(k-2)(t-1)+\sum_{i=1}^t n_i
    \end{align*}
    as desired.
\end{proof}

\section{An observation on edge-heavy independent sets}\label{app:edge_heavy}

The following lemma states that any connected graph contains an independent set incident to many edges.
For a graph~$H$ and a set $S\subseteq V(H)$, we write~$N(S)$ for the set of vertices in $V(H)\setminus S$ which are incident to some vertex in~$S$.
We also write $N[S]=N(S)\cup S$.

\begin{lem}\label{lem:richIndepSet}
    A connected graph $H$ always contains an independent set $I$ which is incident to at least $\max\left\{n-1,2m/\chi(H)\right\}$ edges of $H.$
\end{lem}

\begin{proof}
    For the bound $n-1$, such an independent set can be constructed algorithmically. Start with $I=\{v\}$.
    Iteratively, include a vertex in $V(H) \setminus N[I]$ that has a neighbour in $N(I).$
    Since~$H$ is connected, this process only stops once $N[I]=V(H)$ and thus $H[I,N(I)]$ induces a connected spanning subgraph of $H$, implying $I$ is incident to at least $n-1$ edges.

    For the bound $2m/\chi(H)$, one can take the colour class in an arbitrary proper $\chi(H)$-colouring of~$H$ for which the sum of degrees of its vertices attains the maximum.
\end{proof}

\end{document}